\documentclass[11pt]{amsart}

\usepackage{amssymb,geometry,stmaryrd,color}

\usepackage{amsmath}
\usepackage{amsfonts}
\usepackage{fullpage}

\usepackage{hyperref}
\hypersetup{colorlinks,linkcolor={red},citecolor={blue},urlcolor={blue}}

\usepackage{graphicx}
\usepackage{caption}

\begin{document}

\newtheorem{theorem}{Theorem}[section]
\newtheorem{example}{Example}

\numberwithin{equation}{section}

\def\s{{\bf s}} 
\def\t{{\bf t}} 
\def\u{{\bf u}} 
\def\x{{\bf x}} 
\def\y{{\bf y}} 
\def\z{{\bf z}} 
\def\B{{\bf B}} 
\def\C{{\bf C}} 
\def\K{{\bf K}}
\def\M{{\bf M}}
\def\Nn{{\bf N}}
\def\G{{\bf \Gamma}} 
\def\W{{\bf W}}
\def\X{{\bf X}}
\def\U{{\bf U}}
\def\V{{\bf V}}
\def\Un{{\bf 1}}
\def\Y{{\bf Y}}
\def\Z{{\bf Z}}
\def\P{{\bf P}}
\def\Q{{\bf Q}}
\def\L{{\bf L}}

\def\cB{{\mathcal{B}}} 
\def\cC{{\mathcal{C}}} 
\def\cD{{\mathcal{D}}} 
\def\cG{{\mathcal{G}}} 
\def\cK{{\mathcal{K}}} 
\def\cL{{\mathcal{L}}} 
\def\cR{{\mathcal{R}}} 
\def\cS{{\mathcal{S}}}
\def\cU{{\mathcal{U}}}
\def\cV{{\mathcal{V}}} 
\def\cX{{\mathcal X}}
\def\cY{{\mathcal Y}}
\def\cZ{{\mathcal Z}}

\def\Ea{E_\a}
\def\eps{{\varepsilon}} 
\def\esp{{\mathbb{E}}} 
\def\Ga{{\Gamma}}

\def\lacc{\left\{}
\def\lcr{\left[}
\def\lpa{\left(}
\def\lva{\left|}
\def\racc{\right\}}
\def\rpa{\right)}
\def\rcr{\right]}
\def\rva{\right|}

\def\prst{{\leq_{st}}}
\def\prost{{\prec_{st}}}
\def\prcvx{{\prec_{cx}}}
\def\Rr{{\bf R}}

\def\CC{{\mathbb{C}}}
\def\EE{{\mathbb{E}}}
\def\NN{{\mathbb{N}}} 
\def\QQ{{\mathbb{Q}}} 
\def\PP{{\mathbb{P}}}
\def\ZZ{{\mathbb{Z}}}
\def\RR{{\mathbb{R}}}

\def\Tt{{\bf \Theta}}
\def\Ttt{{\tilde \Tt}}

\def\a{\alpha}
\def\A{{\bf A}}
\def\AA{{\mathcal A}}
\def\hAA{{\hat \AA}}
\def\hL{{\hat L}}
\def\hT{{\hat T}}

\def\claw{\stackrel{(d)}{\longrightarrow}}
\def\elaw{\stackrel{(d)}{=}}
\def\pslaw{\stackrel{a.s.}{\longrightarrow}}
\def\qed{\hfill$\square$}

\newcommand*\pFqskip{8mu}
\catcode`,\active
\newcommand*\pFq{\begingroup
        \catcode`\,\active
        \def ,{\mskip\pFqskip\relax}%
        \dopFq
}
\catcode`\,12
\def\dopFq#1#2#3#4#5{%
        {}_{#1}F_{#2}\biggl[\genfrac..{0pt}{}{#3}{#4};#5\biggr]%
        \endgroup
}

\def\ii{{\rm i}}
\def\S{\mathbf{S}}
\def\F{\mathbf{T}}
\def\W{\mathbf{W}}

\title[log-concave]{An example of a non-log-concave distribution where the difference has a log-concave density}

\author[M.~Wang]{Min Wang}

\address{School of Mathematics and Statistics, Wuhan University, Wuhan 430072, China }

\email{minwang@whu.edu.cn}



\begin{abstract} By the Pr\'ekopa–Leindler inequality, the difference $X-X'$ has a log-concave density provided that $X$ has a log-concave density and $X, X'$ are independent and identically distributed. We prove that the opposite direction does not always hold true by giving an explicit example. \end{abstract}

\maketitle

Log-concave distributions are those for which the density has a concave logarithm. These distributions play an increasingly important role in probability, statistics, optimization theory and other areas of applied mathematics, see the review of Saumard and Wellner \cite{SauWel2014}. They are also necessary for a number of algorithms, e.g. adaptive rejection sampling. Ball, Barthe and Naor \cite{BallBartheNaor2003} investigated log-concave distributions, and tried to understand them through the Fisher information, since if the derivative of $f$ decays at $\infty$, the Fisher information is equal to $\int_\RR f(-\log f)''$, the log-concave distributions are precisely those for which the integrand is pointwise non-negetive. Grechuk, Molyboha and Zabarankin \cite{GrechukMolybohaZabarankin2009} proved that every distribution with log-concave density is a maximum entropy probability distribution with specified mean and Deviation risk measure. Furthermore, log-concave densities are imporant in functional inequalities, for example, Bakry, Barthe, Cattiaux and Guillin \cite{BakryBartheCatGui2008} proved that any log-concave probability measure satisfies both a $L^2$ and a $L^1$ Poincar\'e inequlity.  Ledoux \cite{Ledoux2004} showed in the log-concave situation, $L^2$ and $L^1$ Poincar\'e inequalities are equivalent. Many common probability distributions are log-concave, for example, the normal distribution, the exponential distribution, the logistic distribution, the Laplace distribution, the chi distribution and so on. 

Let $X$ and $X'$ be two independent and identically distributed random variables with log-concave density. It is well-known that, by the Pr\'ekopa–Leindler inequality, log-concavity is preserved by independent summation of log-concave distributed random variables. In particular, $X - X'$ has a log-concave density provided that $X$ has a log-concave density. It is a natural question whether the converse is also valid. We found neither references nor counterexamples in the literature for this kind of reverse to Pr\'ekopa's theorem. Some similar formulations are well known: given independent normally distributed random variables, their sum is normally distributed as well, and Cr\'amer's decomposition theorem says that the converse is also true; given independent random variables of Poisson distribution, their sum has a Poisson distribution as well, and Raikov's theorem says that the converse is also true. It was also proved by Linnik \cite{Linnik1977} that a convolution of normal distribution and Poisson's distribution possesses a similar property.

Simon \cite[Theorem on page 2589]{Simon2011} proved that for standard positive $\alpha$-stable random variable $Z_\a$ with $\a \in (0,1)$, one has the following equivalence
\begin{equation*}
    \log Z_\a \; \text{is log-concave} \; \Longleftrightarrow \; \log Z_\a - \log Z'_\a \; \text{is log-concave} \; \Longleftrightarrow \; \alpha \in (0, 1/2].
\end{equation*}
Then Simon and Barthe \cite[Remark (c) on page 2593]{Simon2011} both wondered whether the log-concavity of $ X - X'$ implies the log-concavity of $X$ in full generality. In the present note we show that it is not true, and we give below an example of a non-log-concave distribution where the difference has a log-concave density. 

\begin{theorem}
    \label{main}
 Let $X_1, X_2, X_3, X_4$ be independent standard normal distributions. Then $X_1 X_2 - X_3 X_4$ has a log-concave density, while $X_1 X_2$ has a non-log-concave density. 
 \end{theorem}

\proof
We first prove that the difference
$X_1 X_2 - X_3 X_4$ has a $\text{Laplace} (0,1)$ distribution, i.e. its probability density is $\frac{1}{2}\exp{(-|x|)}$, by comparing the moment generating function (MGF). The MGF of $X_1X_2$ is, for $|t| < 1$, 
\begin{eqnarray*}
    \EE \lpa e^{tX_1X_2} \rpa &=& \int_{\RR^2} e^{txy} \frac{1}{2\pi} e^{-(x^2 + y^2)/2} dxdy  \\
    &=& \frac{1}{2\pi} \int_0^{2\pi}\int_0^\infty  e^{tr^2\cos \theta \sin \theta} e^{-r^2 /2} rdrd\theta    \\
    &=&  \frac{1}{2\pi} \int_0^{2\pi}\int_0^\infty  e^{r^2(t \cos \theta \sin \theta - 1/2)}  rdrd\theta    \\
    &=&    \frac{1}{2\pi} \int_0^{2\pi} \frac{1}{1-t \sin (2\theta)}  d\theta  = \frac{1}{\sqrt{1-t^2}}.
\end{eqnarray*}
Hence the MGF of $X_1X_2 - X_3X_4$ is $\frac{1}{1-t^2}$ for $|t| < 1$, which is the standard Laplace MGF. It is obvious that the $\text{Laplace} (0,1)$ distribution has a log-concave density. 
The random variables $X_1 X_2$ and $X_3 X_4$ are independent and identically distributed, whose probability density is 
\begin{eqnarray*}
    f_{X_1 X_2} (x) & = & \int_\RR f_{X_1}(\frac{x}{y}) f_{X_2}(y) \frac{dy}{y} = \frac{1}{2\pi}\int_0^\infty e^{-t - x^2/4t} \frac{dt}{t} \\
    &=& \frac{1}{\pi}\int_0^\infty e^{-|x| \cosh(t)}dt =  \frac{1}{\pi} K_0(|x|),
\end{eqnarray*}
where $K_0(x)$ is the modified second class Bessel function, also called the MacDonald or Hankel function.
By \cite[Theorem 2 (b)]{BarPonVuo2011},
the function $K'_\nu(x)/K_\nu(x)$ is strictly increasing on $(0,\infty)$ for all $\nu \in \RR$. Therefore, $K_0(x)$ is a log-convex function on $\RR_+$, thus the random variable $X_1 X_2$ does not have a log-concave density.
\endproof

\bigskip

\noindent
\textbf{Acknowledgements.}
I would like to thank Prof. Thomas Simon and Prof. Fuqing Gao for their encouragements and comments.


\bibliographystyle{plain} 
\bibliography{log-concave}

@article{BallBartheNaor2003,
author = {Keith Ball and Franck Barthe and Assaf Naor},
title = {{Entropy jumps in the presence of a spectral gap}},
volume = {119},
journal = {Duke Mathematical Journal},
number = {1},
publisher = {Duke University Press},
pages = {41 -- 63},
year = {2003},
doi = {10.1215/S0012-7094-03-11912-2},
URL = {https://doi.org/10.1215/S0012-7094-03-11912-2}
}

@article{BarPonVuo2011,
title = {Functional inequalities for modified Bessel functions},
journal = {Expositiones Mathematicae},
volume = {29},
number = {4},
pages = {399-414},
year = {2011},
issn = {0723-0869},
doi = {https://doi.org/10.1016/j.exmath.2011.07.001},
url = {https://www.sciencedirect.com/science/article/pii/S0723086911000442},
author = {Árpád Baricz and Saminathan Ponnusamy and Matti Vuorinen}
}

@article{BakryBartheCatGui2008,
author = {Dominique Bakry and Franck Barthe and Patrick Cattiaux and Arnaud Guillin},
title = {{A simple proof of the Poincaré inequality for a large class of probability measures}},
volume = {13},
journal = {Electronic Communications in Probability},
number = {none},
publisher = {Institute of Mathematical Statistics and Bernoulli Society},
pages = {60 -- 66},
year = {2008},
doi = {10.1214/ECP.v13-1352},
URL = {https://doi.org/10.1214/ECP.v13-1352}
}

@article{GrechukMolybohaZabarankin2009,
 ISSN = {0364765X, 15265471},
 URL = {http://www.jstor.org/stable/40538392},
 author = {Bogdan Grechuk and Anton Molyboha and Michael Zabarankin},
 journal = {Mathematics of Operations Research},
 number = {2},
 pages = {445--467},
 publisher = {INFORMS},
 title = {Maximum Entropy Principle with General Deviation Measures},
 urldate = {2024-03-17},
 volume = {34},
 year = {2009}
}

@incollection {Ledoux2004,
    AUTHOR = {Ledoux, Michel},
     TITLE = {Spectral gap, logarithmic {S}obolev constant, and geometric
              bounds},
 BOOKTITLE = {Surveys in differential geometry. {V}ol. {IX}},
    SERIES = {Surv. Differ. Geom.},
    VOLUME = {9},
     PAGES = {219--240},
 PUBLISHER = {Int. Press, Somerville, MA},
      YEAR = {2004},
   MRCLASS = {58J65 (53C21)},
  MRNUMBER = {2195409},
MRREVIEWER = {Dario Cordero-Erausquin},
       DOI = {10.4310/SDG.2004.v9.n1.a6},
       URL = {https://doi.org/10.4310/SDG.2004.v9.n1.a6},
}

@article{linnik1977,
  title={Decomposition of random variables and vectors},
  author={Linnik, Yu.V. and Ostrovskii, I.V.},
  journal={Translations of Mathematical Monographs},
  year={1977}
}

@article{SauWel2014,
author = {Adrien Saumard and Jon A. Wellner},
title = {{Log-concavity and strong log-concavity: A review}},
volume = {8},
journal = {Statistics Surveys},
number = {none},
publisher = {Amer. Statist. Assoc., the Bernoulli Soc., the Inst. Math. Statist., and the Statist. Soc. Canada},
pages = {45 -- 114},
year = {2014},
doi = {10.1214/14-SS107},
URL = {https://doi.org/10.1214/14-SS107}
}

@article {Simon2011,
    AUTHOR = {Simon, Thomas},
     TITLE = {Multiplicative strong unimodality for positive stable laws},
   JOURNAL = {Proc. Amer. Math. Soc.},
  FJOURNAL = {Proceedings of the American Mathematical Society},
    VOLUME = {139},
      YEAR = {2011},
    NUMBER = {7},
     PAGES = {2587--2595},
      ISSN = {0002-9939,1088-6826},
   MRCLASS = {60E07 (60E05 62E10)},
  MRNUMBER = {2784828},
MRREVIEWER = {Slobodanka\ Jankovi\'{c}},
       DOI = {10.1090/S0002-9939-2010-10697-4},
       URL = {https://doi.org/10.1090/S0002-9939-2010-10697-4},
}
\end{document}